\documentclass[11pt,twoside,openbib]{article}
\usepackage[francais]{babel}
\usepackage{graphicx}
\usepackage{amsmath}
\usepackage{amsfonts}
\usepackage{amssymb}
\usepackage{enumerate}
\newtheorem{theoreme}{Th\'eor\`eme}[section]
\newtheorem{lemme}[theoreme]{Lemme}

\def\HH{{\text {H}}}

\def\HHH{{\mathcal H}}

\title{Global existence of a heavy ball model with a varying nonnegative friction}
\begin{document}
\maketitle
\date{}
\maketitle \centerline{Hassan MCHEIK {\footnote{Lebanese University,
Faculty of Sciences, Beirut, Lebanon. hassan.mcheik@ul.edu.lb
(Hassan MCHEIK), salloum@ul.edu.lb (Zaynab SALLOUM)}} and {Zaynab
SALLOUM}$^{1}$ }
\bigskip 

\noindent {\textbf{Abstract.}} In this paper, we study a dissipative
dynamical system non linear of second order $
\ddot{x}(t)+\lambda(t)\, \dot{x}(t)+\nabla \Phi (x(t))=0,$ with  the
non-negative friction coefficient $\lambda \in
\mathcal{C}([0,+\infty[)$ and the internal constraint  $\Phi \in
\mathcal{C}^2$.\\
\\
\noindent  \textbf{Key-words:} dissipative dynamical system, local minima,
heavy ball with friction depending on time.\\

\noindent{\bf R\'esum\'e} \vskip 0.5\baselineskip
 \noindent {\bf Existence globale d'un syst\`{e}me
mod\'{e}lisant la boule pesante avec un frottement variable
positif.} Dans cette note, nous allons \'{e}tudier le syst\`eme
dynamique dissipatif non lin\'{e}aire du second ordre $
\ddot{x}(t)+\lambda(t)\, \dot{x}(t)+\nabla \Phi (x(t))=0,$ o\`{u}
$\lambda \in \mathcal{C}([0,+\infty[)$ est un c{\oe}fficient du
frottement
positif et $\Phi \in \mathcal{C}^2$ est une contrainte interne.\\

\noindent\textbf{ Mots cl\'{e}s:}  syst\`eme dynamique dissipatif,
point critique,
 fonction convexe, \'{e}quation de la boule pesante d\'{e}pendant du
 temps.

\section{Introduction}
\noindent Dans ce travail, nous \'{e}tudions un  syst\`eme dynamique
dissipatif non lin\'{e}aire, que nous allons nommer HBFT (heavy ball
with friction depending on time),
$$ \ddot{x}(t)+\lambda(t)\, \dot{x}(t)+\nabla \Phi (x(t))=0,$$
 dans le cas  o\`{u} la contrainte externe du frottement $\lambda(t)$ d\'{e}pend du temps avec des conditions initiales sur la
 trajectoire. Notre objectif est de prouver  un r\'{e}sultat
 d'unicit\'{e} et d'existence de la trajectoire avec quelques r\'{e}sultats pr\'{e}liminaires pour
 pr\'{e}parer \`{a} d\'{e}montrer la nature de la convergence (forte ou faible) de la
 trajectoire avec des contraintes externes et internes.\\

 Pour cela, nous commen\c{c}ons par transformer notre probl\`{e}me
 \`{a} un syst\`{e}me du premier ordre dans un espace de Hilbert,
 r\'{e}solu par la m\'{e}thode de Cauchy-Lipschitz.\\

 Le travail est organis\'{e} comme suit. Apr\`{e}s la
 probl\'{e}matique physique
 dans la section 2. et la mod\'{e}lisation math\'{e}matique dans la section 3, nous
 pr\'{e}sentons, dans la section 4,
 nos notations ainsi que nos espaces. Enfin, dans les sections 5 et
 6 nous citons les hypoth\`{e}ses sur les contraintes et nous
 prouvons les th\'{e}or\`{e}mes principaux.
\section{Probl\'{e}matique physique}
Nous consid\'{e}rons une boule pesante avec une contrainte externe
du frottement, qui d\'{e}pend du temps.
 \setcounter{equation}{0}
\noindent Soit $M=(x(t),\Phi (x(t)))$ un point mat\'eriel de masse
$m$, $\Sigma=graphe(\Phi)$ et $\vec{r(t)}$ est  le vecteur de
position \`a l'instant t, d\'{e}fini par  $\vec{r(t)}=(x(t),\Phi
(x(t)).$\\

 Le principe
fondamental de la dynamique (PFD) nous donne:
\begin{equation}\label{0.3.1}
m\ddot{\vec{r}}(t)=\vec{G}+\vec{F}+\vec{R},
\end{equation}
avec   ${\vec{\dot r(t)}}=\big(\dot{x}(t), \nabla\Phi
(x(t))\dot{x}(t)\big)$  est le vecteur vitesse,
$\ddot{\vec{r(t)}}=\big(\ddot{x}(t), \dot{x}(t)
(H_\Phi(x(t)))\ddot{x}(t) +\nabla\Phi (x(t))\ddot{x}(t)\big)$ est le
vecteur d'acc\'{e}leration, o\`{u} $H_\Phi$ est la matrice hessienne
de $\Phi$ en x. $\vec{R}= R\vec{n}$ d\'{e}signe le vecteur de la
r\'eaction sur $\Sigma$, $\vec{F}(t)=-\lambda(t){\dot{\vec{r}}}(t)$
\text{ est le vecteur de la force de frottement}, $\lambda$ est
\text{une contrainte externe de frottement,}
$\vec{G}=\left(0,-mg\right)$ est \text{ la force
 gravitationnelle} et
 $\vec{n}=\frac{1}{\sqrt{1+\left|\nabla\Phi(x)\right|^{2}}}\left(-\nabla\Phi(x),1\right)$ est  \text{ le vecteur
 normal}.\\

Dor\'{e}navant, on s'int\'{e}resse aux cas  dans lesquels les
fonctions $\left|\nabla\Phi(x)\right|$ et $\dot{x}H_{\Phi}\dot{x}$
sont n\'{e}gligeables par rapport \`a 1 et $g$ respectivement
(\textit{i.e.} $\left|\nabla\Phi(x)\right|<<< 1$ et
$\dot{x}H_{\Phi}\dot{x}<<<g$). En effet, le fait que
$\left|\nabla\Phi(x)\right|<<<1$   va nous aider \`{a} trouver
l'ensemble des points limites dans un ensemble bien
d\'{e}termin\'{e} d\'{e}pendant de $\Phi$ (en g\'en\'{e}ral, c'est
l'ensemble des points critiques de $\Phi$) celui-ci est d\'{e}fini
comme l'ensemble des \'el\'ements de $\HHH$ tels que
$\nabla\Phi(x)=0 $. En plus, $\dot{x}H_{\Phi}(x)\dot{x}<<<1$
signifie qu'\`a partir d'un certain temps, la vitesse de la
trajectoire devient tr\`es proche de z\'ero.\\


En faisant les projections canoniques de l'\'{e}quation du principe
fondamental de la dynamique, on obtient les \'equations
diff\'erentielles suivantes:
 \begin{equation}\label{eqt1}
 m\ddot{x}=-\lambda
 \dot{x}-\frac{R}{\sqrt{1+\left|\nabla\Phi(x)\right|^{2}}}\nabla\Phi(x).
 \end{equation}

 \begin{equation}\label{eqt2}
 m(\dot{x}H_{\Phi}\dot{x}+\nabla\Phi(x)\ddot{x})=-mg-\lambda\nabla\Phi(x)\dot{x}+\frac{R}{\sqrt{1+\left|\nabla\Phi(x)\right|^{2}}}.
\end{equation}
\noindent On prend $R>0$, car il  y a toujours des contacts entre la
particule $M(x, \Phi(x))$ et la surface $\Sigma$.\\

\noindent Puis, en projettant suivant la direction normale, on
trouve
\begin{equation}\label{eqt3}
R=\frac{m}{\sqrt{1+\left|\nabla\Phi(x)\right|^{2}}}(g+\dot{x}H_\Phi\dot{x})
\end{equation}

\noindent En utilisant les \'{e}quations syst\'{e}matiques
(\ref{eqt1}) et (\ref{eqt3}), on obtient
\begin{equation}\label{eqt4}
 m\ddot{x}+\lambda
 \dot{x}+\frac{m}{ {1+\left|\nabla\Phi(x)\right|^{2}}}(g+\dot{x}H_\Phi\dot{x})\nabla\Phi(x)=0.
 \end{equation}
\noindent En faisant les aproximations
$\left|\nabla\Phi(x)\right|<<< 1$ et $\dot{x}H_{\Phi}\dot{x}<<<g$,
on obtient $$m\ddot{x}+\lambda \dot{x}+mg\nabla\Phi(x)=0.$$

 \noindent Dans la suite, on s'int\'eresse \`{a}
l'\'equation suivante:
 $$\ddot{x}+\lambda \dot{x}+\nabla\Phi(x)=0.$$
\section{Mod\'{e}lisation Math\'{e}matique}
\setcounter{equation}{0} \noindent On s'int\'{e}resse dans ce papier
au syst\`eme dynamique dissipatif non linéaire $(HBFT)$, qui
mod\'elise l'\'equation de la boule pesante au point mat\'eriel
$M(t)=(x(t),\Phi (x(t))$,
$$ \ddot{x}(t)+\lambda(t) \dot{x}(t)+\nabla \Phi (x(t))=0, \lambda(t) >0,$$
o\`u la contrainte  $\lambda$ est une fonction correspondante \`a un
coefficient de frottement, avec des conditions initiales classiques
et des contr\^{o}les optimaux sur les contraintes internes et
externes $\Phi$ et $\lambda$, respectivement.\\

 Intuitivement,
l'interpretation m\'{e}canique de ce syst\`{e}me d\'{e}pend de
l'\'{e}nergie initial donn\'{e} par
$E_0=\frac{1}{2}\left|\dot{x_0}\right|^{2}+\Phi (x_0)$, o\`{u}
$x_0=x(0)$ et $\dot{x_0}=\dot{x}(0)$ sont respectivement le point de
d\'{e}part et
  la v\'{e}locit\'{e}  au temps initial $t_0=0$. D'autre part, la position du point $M(t))=(x(t),\Phi
(x(t))$ d\'{e}pend de la contrainte  externe $\lambda(\cdot)$ et
tend de proche en proche vers le point $M^*$, qui est en
g\'{e}n\`{e}ral un point stable d\'{e}pend de la caract\'{e}risation
de $\Phi$. Les r\'esultats concernant la convergence de la solution,
par exemple  vers un point critique, sont donn\'{e}s dans plusieurs
cas, incluant le cas o\`u $\Phi$ est une fonction convexe; une
fonction de morse ou une fonction  ayant un nombre fini
d'intervalles de stabilit\'{e}. La trajectoire du syst\`{e}me $(HBFT)$ converge asymptotiquement vers un minimum local de $\Phi$.\\

Dans la suite, on d\'efinit le syst\`eme dynamique (HBFT) suivant
[cf. \ref{Al},\ref{AC},\ref{ACZ},\ref{ACR}]:
\begin{equation}\label{1.}
\left\{
\begin{array}{c}
\  \ddot{x}(t)+\lambda (t)\dot{x}(t)+\nabla \Phi (x(t))=0\text{\ }   \\
\ \lambda :[0,+\infty \lbrack \rightarrow R_{+}\text{ continue. \  }\\
\end{array}
\right.
\end{equation}
De plus, l'\'energie de ce  syst\`eme est d\'{e}finie  par
 $E(t)=\frac{1}{2}\left|\dot{x}(t)\right|^{2}+\Phi (x(t)) $, dont la
d\'eriv\'ee  est  $ \dot{E}(t)=-\lambda
(t)\left|\dot{x}(t)\right|^{2}.$ On remarque que l' \'energie est
d\'ecroissante polynomialement en fonction du temps.

\section{Notations, espaces et normes}
Dans toute la suite, $\mathcal{H}$ d\'{e}note un espace de Hilbert,
muni du produit scalaire $<\cdot,\cdot>$ et de la norme associ\'{e}e
$|\cdot|$.  Nous utilisons l'espace $L^p$, qui  est un espace
vectoriel de classes des fonctions dont la puissance d'exposant $p$
est int\'{e}grable au sens de Lebesgue, o\`{u} $p$ est un nombre
r\'{e}el strictement positif. 
On d\'{e}signe par $\text{argmin}\Phi$ l'ensemble non vide des
points critiques de $\Phi$,\textit{ i.e.} $\text{argmin}\Phi=:\{x
\in \mathcal{H}, \nabla \Phi (x(t))=0 \}$.
\subsection{Hypoth\`{e}ses sur les contraintes}
Dans cette partie, nous commen\c{c}ons par pr\'{e}ciser des
hypoth\`{e}ses sur les fonctions $\Phi :\mathcal{H}\rightarrow
\mathbb{R}$ et
$ \lambda :[0,+\infty \lbrack \rightarrow \mathbb{R}_{+}$. \\ 
Les hypoth\`eses standards sur $\Phi$ (voir [\ref{AGR}]) seront
not\'{e}es $(\text{H}_{\Phi})$, et donn\'{ee}s par:
\begin{equation*}
(\text{H}_{\Phi })\left\{
\begin{array}{l}
\text{H}_{\Phi }^{\tiny 1}.\quad \text{ la fonction  } \ \Phi \ \text{est de classe} \ C^{1}  \ \text{sur} \ \mathcal{H}.\\
\text{H}_{\Phi }^{\tiny 2}.\quad \text{ la fonction} \ \nabla \Phi \ \text{est Lipschitzienne sur les parties born\'ees de} \ \mathcal{H}.\\
\text{H}_{\Phi }^{\tiny 3}.\quad \text{ la fonction} \ \Phi \
\text{est born\'ee inf\'erieurement.}
\end{array}
\right.
\end{equation*}
%

\noindent Dans ce papier, nous proposons des conditions sur la
contrainte externe $\lambda$ qui d\'{e}pend du temps. Ces
hypoth\`eses sont not\'{e}es $(\text{H}_{\lambda })$  et donn\'{e}es
par :
\begin{equation*}
\text (H_{\lambda })\left\{
\begin{array}{l}
\text{H}_{\lambda }^{\tiny 1}.\quad \lambda \text{ est une
fonction continue sur $[0,\infty[$.}\\
\text{H}_{\lambda }^{\tiny 2}.\quad\text{Il existe }  M>0 \text{ et
} t_{1}\geq 0, \text{ tel que  pour tout }  \, t\geq t_{1},
\text{on a } \lambda(t)\leq M,\\
 \qquad{\textit{i.e. } } \underset{t\rightarrow \infty }{\limsup \text{ }}
 \lambda(t)<\infty.\\
\end{array}
\right.
\end{equation*}
%
\subsection{R\'{e}sultats principaux.}
Dans cette section, nous  commen\c{c}ons  par  \'{e}tablir   les
r\'{e}sultats d'existence  et d'unicit\'{e} de solutions du
syst\`{e}me (HBFT), en utilisant le th\'{e}or\`{e}me de
Cauchy-Lipshitz et par application du th\'{e}or\`{e}me d'Ascoli (cf.
[\ref{AGR}]). En effet,  H. Attouch, X. Goudou et P. Redont
[\ref{AGR}] ont \'{e}tudi\'{e} le syst\`{e}me (HBF), o\`{u} la
contrainte $\lambda$ est une constante, et $\Phi$ v\'{e}rifie  les
conditions standards $(\text{H}_{\Phi})$.

Nous commen\c{c}ons par annoncer le premier r\'{e}sultat:

\begin{theoreme}\label{th1.1}
Soient $\Phi: \mathcal{H}\rightarrow \mathbb{R}$ et $ \lambda
:[0,+\infty \lbrack \rightarrow \mathbb{R}_{+}$ deux  fonctions
v\'erifiant les hypoth\`{e}ses $(\text{H}_{\Phi })$  et
$\text{H}_{\lambda}^1$, alors on a les propri\'et\'es suivantes :
\begin{enumerate}[(i)]
\item Pour tout  $(x_{0},\dot x_{0} )$ dans $\mathcal{H}\times \mathcal{H}$,
 il existe une solution unique $x$ de classe $C^{2}$ de syst\`{e}me $(\ref{1.})$ d\'efinie sur l'intervalle $[0,+\infty[
 $ et v\'erifiant les conditions initiales $x(0)=x_{0},  \dot{x}(0)=\dot x_{0}.$
\item Pour toute trajectoire $x$ de $(\ref{1.})$, l'\'energie correspondante $E$ est d\'ecroissante polynomialement sur $[0,+\infty $ $[$ \ et converge vers $E_{\infty }$ et de plus:
$$
 \dot{x}\in L^{\infty }([0,\infty \lbrack ;\HHH)\ \text{et }\sqrt{\lambda}\dot{x}\in L^{2}([0,\infty \lbrack ;\HHH).
$$
\end{enumerate}
\end{theoreme}
%
%
%

\subsection*{Preuve du Th\'eor\`eme \ref{th1.1}}
\noindent
 \begin{enumerate}[(i)]
\item Pour toute condition initiale $(x_{0},\dot{x}_{0})\in \HHH\times \HHH$, d'apr\`es le th\'eor\`eme de Cauchy-Lipschitz plus les hypoth\`eses $\HH_{\Phi}$ et
$\HH_{\lambda}^1,$ le syst\`eme   $(\ref{1.})$ poss\`ede une
solution unique de classe $C^{2}$, d\'efinie sur l'intervalle
$[0,T_{\max }[$, avec $0<T_{\max }\leq \infty. $ Dans le but de
trouver $T_{\max }=\infty $, montrons que $ \dot{x}$ est born\'ee.\\
L'\'equation $(\ref{1.})_1$ et la r\'egularit\'e de $\Phi$ et $\HH_{\lambda}^1$ impliquent que $x$ est de classe $C^{2}$ sur $[0,T_{\max }[$.\\
En effet, on transforme le syst\`eme $(HBFT)$ en un syst\`eme
r\'eduit dans $\HHH\times \HHH$ de la forme :
\begin{equation}\label{2.1} \dot{Y}(t)=F(Y(t),t),\qquad Y_{0}=\left(
\begin{array}{l}
x_{0} \\
\dot{x}_{0} %
\end{array}%
\right) \in  \ \HHH \times \HHH,
\end{equation}
avec
$$
Y(t)=\left(
\begin{array}{l}
x(t) \\
 \dot{x}(t)%
\end{array}%
\right) \ \ et\ \ \ \ F(u,v,t)=\left(
\begin{array}{l}
v \\
-\lambda (t)v-\nabla \Phi (u)%
\end{array}%
\right).
$$
Puisque $\nabla \Phi $ est localement Lipschitzienne et $\lambda$
est continue, alors le th\'eor\`eme de Cauchy-Lipschitz nous affirme
l'existence et l'unicit\'{e} d'une solution locale de  probl\`eme
(\ref{2.1}).
Soit $x$ la solution maximale du syst\`eme $(HBFT)$, d\'efinie sur
$[0,T_{\max }[$, avec $0<T_{\max }\leq \infty .$ Tout d'abord,
l'\'energie de ce syst\`{e}me dissipatif  est donn\'{e} par:
\begin{equation*}
 E(t)=\frac{1}{2}\left|\dot{x}(t)\right|^{2}+\Phi (x(t)).
 \end{equation*}
En d\'erivant l'\'{e}quation de l'\'{e}nergie $E$ puis en utilisant
$(\ref{1.})$, on obtient:
\begin{equation}\label{2.2}
\dot{E}(t)=-\lambda (t)\left|\dot{x}(t)\right|^{2}.
\end{equation}
Puisque  $\lambda (t)\geq 0$ pour tout $  t\geq 0$, alors la
fonction $E$ est d\'ecroissante sur $[0,T_{\max }[$, et pour tout
$t$ appartient \`a  $[0,T_{\max }[$,
 on a  $E(t)\leq E({0}).$ Cette in\'{e}galit\'{e} est \'{e}quivalente \`a 
\begin{equation}\label{2.3}
\frac{1}{2}\left|\dot{x}(t)\right|^{2}+\Phi (x(t))\leq
\frac{1}{2}\left|\dot{x}_{0}\right|^{2}+\Phi (x_{0}).
\end{equation}
De plus, puisque $\Phi $ est minor\'ee (hypoth\`{e}se
$\text{H}_{\Phi }^{\tiny 3}$), alors on peut affirmer que,
$$\left\Vert  \dot{x}\right\Vert _{\infty }=\underset{[0,T_{\max }[}{\sup }%
\left\vert  \dot{x}(t)\right\vert <\infty.$$
Montrons maintenant que $T_{\max}=\infty .$\\

\noindent En effet, supposons que $T_{\max }<\infty$, alors
d'apr\`{e}s le th\'{e}or\`{e}me des accroissements finis:
 $$
\forall t,\, t' \in[0,T_{\max }[, \left\vert x(t)-x(t^{\prime
})\right\vert \leq \left\Vert  \dot{x}\right\Vert _{\infty
}\left\vert t-t^{\prime }\right\vert.
 $$

\noindent D'autre part, puisque $T_{\max }<\infty$ et $x$ est
continue sur $[0,T_{\max }[$, alors $\underset{t\rightarrow T_{\max
}}{\lim }x(t)=x_{\infty }$ existe, alors $x$ et $ \dot{x}$ sont
born\'ees sur $[0,T_{\max }[$, et par l'\'equation $(\ref{1.})$, la
fonction $\ddot{x}$ est born\'ee sur cet intervalle. Alors
$\underset{t\rightarrow T_{\max }}{\lim } \dot{x}(t)=\dot{x}_{\infty
}$ existe. Par suite, on note  $\underset{t\rightarrow
T_{\max }}{\lim }(x(t),\dot{x}(t))=(x_{\infty},\dot{x}_{\infty })$.\\
De nouveau, on  applique le th\'eor\`eme de Cauchy-Lipschitz avec
les  conditions initiales obtenues ($x_{\infty},\dot{x}_{\infty })$.
Dans ce cas, $T_{\max }$ est le  point de d\'{e}part et la solution
maximale est d\'efinie sur un intervalle $[ T_{\max }, T'_{\max }[$.
Par suite, on   trouve une  solution locale de $(\ref{1.})$,
d\'{e}finie  sur un intervalle $[0,T'_{\max }[$ plus large que
$[0,T_{\max }[$ avec $T_{\max }< T'_{\max }$, ce  qui donne une
contradiction, donc $T_{\max }=\infty .$
\item On a prouv\'e que $E$ est d\'ecroissante sur $[0,\infty[$
 et d'apr\'es $\HH_{\Phi}$, $\Phi$ est minor\'ee et de plus $E(t)\geq \Phi (x(t))$, alors l'\'{e}nergie  $E$ est
  minor\'ee. Par cons\'equence  $\underset{t\rightarrow \infty }{\lim }E(t)=E_{\infty }$ existe,
  avec $E_{\infty }\in\mathbb{R}.$\\
En utilisant $(\ref{2.3})$ et l'hypoth\`{e}se $\HH_{\Phi}^3$, on
obtient que:
\begin{equation}\label{2.4}
\frac{1}{2} \left|\dot{x}(t)\right|^{2}\leq \frac{1}{2}%
\left|\dot{x}_{0}\right|^{2}+\Phi (x_{0})-\inf \Phi, \textup{ pour
tout }   t\geq  0,
\end{equation}
alors,
$$ \dot{x}\in L^{\infty }([0,\infty \lbrack ;\HHH).$$
De $(\ref{2.2}),$ on a pour tout $0\leq t< \infty$,
 $$
\int_{0}^{t}\lambda (s)\left\vert  \dot{x}(s)\right\vert ^{2}ds=-\int_{0}^{t}\dot{E}(s)ds.
$$
Par suite, \ $$\int_{0}^{\infty }\lambda (s)\left\vert  \dot{x}(s)\right\vert^{2}ds= -\int_{0}^{\infty }\dot{E}(s)ds
 =-\Big(E_{\infty }-E_{0}\Big).$$
Donc
$$
  \sqrt{\lambda(.)}\dot{x}\in L^{2}([0,\infty \lbrack ;\HHH).
 $$
\end{enumerate}
  \begin{lemme}\label{Lemme1.1}
Soit $f: [0;\infty[\rightarrow \HHH$ une fonction de classe
$C^{1}([0;\infty[;\HHH)$  et v\'erifiant $ f\in L^{2}([0,\infty
\lbrack ;\HHH)\cap L^{\infty }([0,\infty\lbrack ;\HHH)$, et $
\dot{f}\in L^{\infty }([0,\infty\lbrack ;\HHH)$, alors on a
$\underset{t\rightarrow \infty }{\lim }f(t)=0.$
\end{lemme}
\subsection*{\textbf{Preuve du Lemme \ref{Lemme1.1}}}
\noindent On d\'{e}finit la fonction $h: [0;\infty[\rightarrow
\mathbb{R}^+$ par $h(t)=\|{f(t)}\|^2$. Il suffit de d\'{e}montrer
 que $\underset{t\rightarrow \infty }{\lim }h(t)=0$.\\
En effet, supposons que   $\underset{t\rightarrow \infty }{\lim }\,h(t)\neq0$.\\

\noindent Tout d'abord, nous allons d\'{e}montrer que
$\underset{t\rightarrow \infty }{\liminf}\, h(t)=0$.  Sinon, c'est
\`{a} dire $\underset{t\rightarrow \infty }{\liminf}\, h(t)>0,$ par
suite il existe $\ \varepsilon >0$ et $t_0>0$ tels que
  pour tout $t \geq t_0$, on a $h(t) \geq \varepsilon$. Alors,
  $$\int_{0}^{\infty}\left\|f(t)\right\|^{2}dt=\int_{0}^{\infty}h(t)dt\geq\int_{t_0}^{\infty}h(t)dt\geq \int_{t_0}^{\infty} \varepsilon dt=\infty, $$
ce qui contredit avec le fait que $f\in L^{2}([0,\infty \lbrack
;\HHH)$.\\

\noindent Par suite, $0=\underset{t\rightarrow \infty }{\liminf}\,
h(t)< \underset{t\rightarrow \infty }{\limsup}\, h(t)$. Donc, il
existe $ \varepsilon>0$, tel que $\underset{t\rightarrow \infty
}{\limsup}\, h(t)>\varepsilon$. En utilisant la continuit\'{e} de la
fonction $h$ et le th\'{e}or\`{e}me des valeurs interm\'{e}diaires,
ils existent deux suites $(s_n)_n$ et    $(t_n)_n$, qui tendent vers
l'infini, telles que $h(t_n)=\varepsilon$ et
$h(s_n)=\frac{\varepsilon}{2}$. \\

\noindent De  plus, $h$ est d\'{e}rivable et
$\dot{h}(t)=2<f(t),\dot{f}(t)>$. Ensuite, $\dot{h}\in L^{\infty
}([0,\infty\lbrack) $ car $|\dot{h}(t)|\leq \|\dot{f}\|_\infty
\|{f}\|_\infty<\infty.$

\noindent En cons\'{e}quence, $h$ et $c-$lipshtitizienne et
$\frac{\varepsilon}{2}=|h(t_n)-h(s_n)|\leq c|t_n-s_n|.$ Donc il
existe    $\eta >0$ telle que $|t_n-s_n|> \eta$.\\

\noindent Soit $\tau_n=\displaystyle\frac{s_n+t_n}{2}$, alors
$\tau_n$ converge vers l'infini. On d\'{e}duit qu'il existe
$\delta=\frac{\eta}{4}
>0$ tel que pour tout  $ t$, avec $\ \left\vert \tau_{n}-t\text{ }\right\vert
\leq \delta$, on a $h(t)\geq \frac{\varepsilon}{2}. $\\

\noindent Alors $$\int_{0}^{\infty}\left\|f(t)\right\|^{2}dt=
\int_{0}^{\infty}h(t)\,dt\geq\sum_{n\geq1}
 \int_{t_{n}-\delta}^{t_{n}+\delta}h(t)\,dt\geq\sum_{n\geq1}
 \int_{t_{n}-\delta}^{t_{n}+\delta}\frac{\varepsilon}{2}\, dt=\infty,$$ donc $f\notin L^{2}([0,\infty \lbrack ;\HHH),$ ce
 qui contredit avec l'hypoth\`{e}se. Donc, $\underset{t\rightarrow \infty }{\lim
 }h(t)=0.$\\
\begin{theoreme}\label{th2.1}
Sous les conditions du th\'{e}or\`{e}me  (\ref{th1.1}) et  si de
plus $x$ est dans $L^{\infty }([0,\infty \lbrack ;\HHH)$, alors on
a:
\begin{enumerate}
\item[i.] $\ddot{x}\in  L^{\infty }([0,\infty \lbrack ;\HHH)$.
\item[ii.]  Si ${\lambda}$ est d\'{e}rivable et $\dot{\lambda}$ est born\'e, alors
$$\underset{t\rightarrow \infty }{\lim \text{ }} \sqrt{\lambda(t)}\dot{x}(t) =0.$$
\end{enumerate}
\end{theoreme}
\subsection*{Preuve du Th\'eor\`eme \ref{th2.1}}
\noindent Puisque $x$ appartient \`a  $L^{\infty }([0,\infty \lbrack
;\HHH)$, $\dot{x}\in L^{\infty}([0,\infty \lbrack ;\HHH)$ et
$\lambda(.)$ est born\'e, alors d'apr\`es  Th\'{e}or\`{e}me
$(\ref{th1.1})_{ii.}$ on a  $\sqrt{\lambda(.)}\dot{x}\in
L^{\infty}([0,\infty \lbrack ;\HHH)$.\\
En utilisant l'\'equation $(\ref{1.})$, et $\nabla \Phi $ qui est born\'ee
sur les parties born\'ees de $\HHH$, on trouve que $ \ddot{x}\in L^{\infty }([0,\infty \lbrack;\HHH).$\\
 Montrons maintenant que $\underset{t\rightarrow \infty }{\lim
 }\sqrt{\lambda(t)}\dot{x}(t)=0.$  En effet, soit $h(t)=\lambda(t)\left|\dot{x}(t)\right|^{2}$. Il suffit de montrer que
 $ h\in L^{2}([0,\infty \lbrack ;\mathbb{R})\cap L^{\infty}([0,\infty \lbrack ;\mathbb{R})$, et $ \dot{h}\in L^{\infty }([0,\infty
\lbrack ;\mathbb{R})$.\\
 Tout d'abord, on a $ h\in L^{2}([0,\infty \lbrack
;\mathbb{R})$ car $${\int_0^\infty h^2(t) dt =\int_0^\infty
\lambda^2(t)\left|\dot{x}(t)\right|^{4} dt \leq  \| \lambda\|_\infty
\|\dot{x}\|_\infty^2 \int_0^\infty
\lambda(t)\left|\dot{x}(t)\right|^{2} dt< \infty }.$$ De plus,
$$ \dot{h}(t)=\dot{\lambda}(t)\left|\dot{x}(t)\right|^{2}+2\lambda(t)\left\langle\dot{x}(t),
\ddot{x}(t)\right\rangle\in L^{\infty }([0,\infty\lbrack
;\mathbb{R}),$$ car $\lambda$, $\dot{\lambda}$, $\dot{x}$ et
$\ddot{x}$ sont
born\'{e}es.\\
 En utilisant  le Lemme (\ref{Lemme1.1}), on obtient aue
$\underset{t\rightarrow \infty }{\lim }h(t)=0$ et par suite
$\underset{t\rightarrow \infty }{\lim
 }\sqrt{\lambda(t)}\dot{x}(t)=0.$
\\

\noindent \textit{\textbf{Remarque}}. \textit{Nous indiquons que
lorsque la trajectoire $x$ du syst\`{e}me (HBFT) est pr\'{e}compact
pour la norme induite de $\HHH$. Le r\'{e}sultat du  Th\'eor\`eme
\ref{th2.1}, $\underset{t\rightarrow \infty }{\lim \text{ }}
\sqrt{\lambda(t)}\dot{x}(t) =0$ peut \^{e}tre obtenu comme
cons\'{e}quence du principe d'invariance de LaSalle (cf. [\ref{B}],
[\ref{H}]).}

 \end{document}